%% file: A5action.tex
\documentclass{article}
\usepackage{epsfig}
\usepackage{amsmath}
\usepackage{amssymb}
\newcommand{\qbar}{\overline{\mathbb{Q}}}
\newcommand{\cur}{\mathcal{C}}
\newcommand{\belyi}{Bely\u{\i} }

\begin{document}
\centerline{\bf An $A_5$ Galois action on a class}
\centerline{\bf of elliptic dessins d'enfants}

\baselineskip=16pt
\vspace{0.8cm}
\centerline{David S. Burggraf}
\baselineskip=13pt
\centerline{ \em Department of Mathematics, University of British
Columbia}
\baselineskip=12pt
\centerline{\em Vancouver, B.C.\ V6T 1Z2, Canada}
\vspace{0.9cm}

\input{abstract}
\input{introduction}
\input{EllipticDessins}
\input{MonodromyGroups}
\input{GaloisActions}

\bibliographystyle{alpha}
\bibliography{A5action.bib}
\end{document}

%% file: abstract.tex
\begin{abstract}

As noted in \cite{jones1}, only examples of solvable Galois actions on classes of dessins d'enfants have been explicitly constructed. In this paper an action of the nonsolvable alternating group $A_5$, represented as a Galois group, on a class of elliptic dessins is presented.

\end{abstract}

%% file: introduction.tex
\section{Introduction}
A \belyi function $\beta$ defined on a compact Riemann surface $X$ is a
meromorphic covering of the Riemann sphere $\Sigma = \mathbb{C} \cup 
\{\infty\}$ with all branch points contained in $ \{ 0,1,\infty \}$. \belyi functions $\beta$ that satisfy the additional condition that all points in the pre-image of 1 are ramified with ramification order equal to  2, are called {\em clean} \belyi functions. It is
common practice to extend the definition of a \belyi function to a projective
algebraic curve $\cur$ by identifying $\cur$ with its underlying compact
Riemann surface. One can associate a bipartite graph $\mathcal G(\beta)$ to the \belyi pair $(\cur,\beta)$ by using
$\beta$ to lift the closed unit interval $I=[0,1] \subset \mathbb{C} \cup
\{\infty\}$; label the points in $\beta^{-1}(I)$ lying above 0 and 1,
respectively, by black and white vertices, respectively. Each connected
component in the inverse image $\beta^{-1}(0,1)$ is an {\em edge segment} joining two
vertices of opposite type and each connected component of $\cal{C}$$\setminus
\beta^{-1}(I)$ is a simply connected region, which will be referred to as a face. A bipartite map $\mathcal M$ on $\mathcal C$ consists of a bipartite graph $\mathcal G$ lying on $\mathcal C$ together with the faces of $\mathcal C \setminus \mathcal G$. A bipartite map on a compact Riemann
surface obtained by lifting $I$ by a clean \belyi funtion is denoted by $\mathcal D(\beta)$, and belongs to an important class of combinatorial objects known as dessins d'enfants. 

\belyi's theorem \cite{bel} states that a projective algebraic curve $\cur$ is defined over the algebraic numbers $\qbar$ if and only if there exists a \belyi function $\beta:\cur \rightarrow$ $\mathbb{P}^{1}(\mathbb{C})=\mathbb{C}  \cup
\{\infty\}$. Inspired by Bely\u{\i}'s theorem, Alexander Grothendieck \cite{gro} observed that the absolute Galois group $\Gamma=$Gal$(\qbar/\mathbb{Q})$ acts naturally on \belyi pairs $(\cur,\beta)$ and induces a faithful action on the corresponding dessins. This action of $\Gamma$ on dessins lends direct geometrical insight into the intractible group $\Gamma$ and is the basis of an ambitious program started by Grothendieck aiming at a complete description of the absolute Galois group.

In the survey paper \cite{jones1} some explicit examples of nontrivial Galois actions on dessins lying on: the Riemann sphere, elliptic curves, and hyperelliptic curves are given. In these examples, the Galois groups acting were solvable, namely: cyclic groups and the metabelian affine group $AGL_1(p)$. This leads to the natural generalization suggested by the author in the last sentence of the paper: to produce explicit examples of nontrivial Galois actions on classes of dessins where more complicated Galois groups are acting such as nonsolvable groups or solvable groups with unbounded derived length. In this paper, the simplest nonsolvable group $A_5$ is explicitly shown to have a nontrivial Galois action on a class of elliptic dessins.

%% file: EllipticDessins.tex
\section{Elliptic Dessins}
\begin{figure}[htb]
\begin{center}
\begin{picture}(100,210)
\put(0,-65){\epsfig{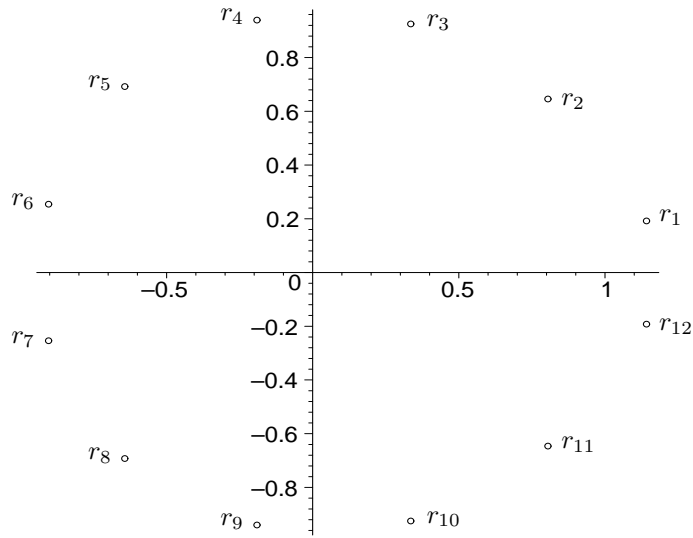}}
\put(170,115){$r_1$}
\put(133,159){$r_2$}
\put(82,189){$r_3$}
\put(4,191){$r_4$}
\put(-46,166){$r_5$}
\put(-75,121){$r_6$}
\put(-75,69){$r_{7}$}
\put(-46,25){$r_{8}$}
\put(4,0){$r_{9}$}
\put(82,1){$r_{10}$}
\put(133,29){$r_{11}$}
\put(170,74){$r_{12}$}
\end{picture}
\end{center}
\caption{ The labelling of the zeros of the polynomial $f(x)$}
\label{fzeros}
\end{figure}

It can be verified using a computer algebra program such as \cite{GAP98} that the separable polynomial $f(x)=x^{12}-\frac{12}{11}x^{11}+1$ is irreducible over $\mathbb{Q}$, and
the Galois group $G = \mbox{Gal}(\mathbb Q,f)$ of its splitting field over $\mathbb{Q}$ is $S_{12}$. The zeros of $f(x)$ are displayed in Figure \ref{fzeros} and labelled for later reference.

\noindent Consider the elliptic curve
\begin{equation}
\label{Eijk}
E_{i,j,k} = \left\{ [x,y,z]\in \mathbb{P}^2(\mathbb{C}) \, | \,
y^2z-(x-r_iz)(x-r_jz)(x-r_kz)=0 \right\},
\end{equation}
where $r_i,r_j$, and $r_k$ are distinct zeros of the polynomial $f(x)$. It is customary, when defining a \belyi function on a projective algebraic curve $\mathcal{C}$, to ``affinize'' the curve $\mathcal{C}$ by considering only the the part of $\mathcal{C}$ contained in the open neighbourhood of $\mathbb{P}^2(\mathbb C)$ where $z \neq 0$. Thus, only the homogeneous coordinates of the form $[x,y,1]$ in $\mathcal C$ are considered, and identifying these with the affine coordinates $(x,y)$, yields the associated affine curve $\mathcal{C}^{a\!f\!f}$. The projective curve is then easily recovered from the affine one by "homogenizing" it's defining polynomial. Therefore,
we will restrict our attention to the affine elliptic curve 
\begin{equation*} 
E_{i,j,k}^{a\!f\!f} =\left\{ (x,y)\in \mathbb{C}^2 \, | \, y^2=(x-r_i)(x-r_j)(x-r_k) \right\}
\end{equation*} 
and make the identification $E_{i,j,k}=E_{i,j,k}^{a\!f\!f} \cup \infty$, where $\infty$ corresponds to the point $[0,1,0]\in \mathbb{P}^2(\mathbb C)$.
The mapping $\pi^{i,j,k} :E_{i,j,k} \rightarrow \Sigma$ defined by $\pi^{i,j,k}(x,y)=x$ is a degree 2 covering of $\Sigma$ with branch points $r_i,r_j,r_k,$ and $\infty$, bearing in mind that $(\pi^{i,j,k})^{-1}(\infty)=[0,1,0]=\infty$. These branch points are in turn sent to $0,0,0,$ and $\infty$, respectively, by the polynomial $f(x)$. In addition, $f(x)$ has two ramification points $0$ and $1$ lying above the branch points $1$ and $\frac{10}{11}$, hence $f\circ\pi^{i,j,k}$ has a total of four branch points $\{0,\frac{10}{11},1,\infty\}$. We will now make use of the well known \belyi function $\beta_{m,n}$ (see \cite{jones2}), where $m=10, n=1$, as it is ideally suited to handle the present situation. The function $\beta_{m,n}:\Sigma \rightarrow \Sigma$ is defined by  
\begin{equation*} 
\beta_{m,n}(x)=\frac{(m+n)^{m+n}}{m^mn^n}x^m(1-x)^n.
\end{equation*}
It can be easily verified that $\beta_{m,n}$ is unramified outside $\{0,\frac{m}{m+n},1,\infty\}$ and sends this set into $\{0,1,\infty\}$. Therefore, the composition $\phi= \beta_{10,1}\circ f \circ \pi^{i,j,k}$ is a \belyi function, and furthermore, $\beta^{i,j,k}= \beta_{1,1}\circ\phi$ is a clean \belyi function, i.e., has ramification order equal to 2 at each point in the pre-image of 1. The lift of the unit interval $I=[0,1]\subset\Sigma$ by $\beta^{i,j,k}$ gives rise to an elliptic dessin $\mathcal D(\beta^{i,j,k})$ (see Figure \ref{ellipticDessin}).

%% file: MonodromyGroups.tex
\section{Monodromy Groups}
Let $\beta:\mathcal C \rightarrow \Sigma$ be a branched covering of the Riemann sphere, and suppose $B=\{b_0,\ldots,b_k\} \subset \Sigma$ is the set of branch points. Fix a base point $p \in \Sigma-B$, and let $\Omega= \beta^{-1}(p)$ denote the fiber above $p$ in $\mathcal C$. The fundamental group $\Pi=\pi_1(\Sigma-B,p)$ acts naturally on $\Omega$ and has a presentation given by
\begin{equation*}
\Pi = \langle \alpha_0,\ldots,\alpha_k \, |\, \prod_{i=0}^k{\alpha_i} =1 \rangle,
\end{equation*}
where each generator $\alpha_i$ is the homotopy class of a suitably chosen loop based at $p$, that winds once around $b_i$. For any element $p_r \in \Omega$, each element $\alpha_i \in \Pi$ lifts to a unique homotopy class of paths in $\mathcal C$, such that any representative path in the homotopy class starts at $p_r$ and ends at some common point $p_s \in \Omega$. This defines a bijective mapping $g_i:\Omega \rightarrow \Omega$, which gives a permutation representation $\rho: \Pi \rightarrow  S_{\Omega}$, defined by $\rho(\alpha_i)=g_i$, where $S_{\Omega}$ is the symmetric group on the set $\Omega$. The permutation $g_i$ ($0 \leq i \leq k$) has $n_i=|\beta^{-1}(b_i)|$ disjoint cycles, each having cycle length $l_{ir}$ ($1 \leq r \leq n_i$). Each disjoint cycle of $g_i$ corresponds to a ramification point of order $l_{ir}$ on $\mathcal C$ lying above $b_i$, where $l_{ir}$ sheets of $\mathcal C$ come together. The resulting permutation group $G(\beta)=\rho(\Pi)$ is called the monodromy group of the branched covering $\beta$. In the case that $\beta$ is a \belyi function, we  have $B \subset \{0,1,\infty\}$, and $G(\beta)$ generated by the two permutations $g_0$ and $g_1$. Moreover, if $\beta$ is a clean \belyi function, all ramification points above 1 have ramification order equal to 2, hence all disjoint cycles of $g_1$ have length $l_{1r}=2$ $ (1 \leq r \leq n_1$), and the monodromy group $G(\beta)$ has the following partial presentation
\begin{equation*}
G(\beta) = \langle g_0,g_1 \, |\, g_0^l=g_1^2 =1, \mbox{ETC} \rangle,
\end{equation*} 
where $l$ is the least common multiple of $\displaystyle \{l_{0r}\}_{r=1}^{n_0}.$

The monodromy group $G(\beta)$ has an important connection to the graph $\mathcal G(\beta)$ associated to $\beta$; in fact, when $\beta$ is a clean \belyi function, the conjugacy class of the monodromy group $G(\beta)$ in $S_{\Omega}$ completely determines the graph-isomorphism class of $\mathcal G(\beta)$ and vice versa. Suppose the base point $p \in \Sigma - B$ is chosen to lie in the open interval $(0,1)$, then each edge segment $e_r$ of $\mathcal G(\beta)$ contains exactly one point $p_r \in \Omega$ and is incident to a unique black vertex at one end and a white vertex at the other. Each disjoint cycle of the monodromy generator $g_0$ describes how the point $p_r$ is cyclically permuted in $\Omega$ around a ramification point $v$ lying above 0, and thus describes how the edge segment $e_r$ is cycled (according to the orientation of $\mathcal C$) around its incident black vertex $v$. The valency of $v$ is then given by the length $l_{0r}$ of the disjoint cycle. This describes the connection between the black vertices of $\mathcal G(\beta)$ and the disjoint cycles of $g_0$, and a similar connection exists between the disjoint cycles of $g_1$ and the white vertices. In the case that $\beta$ is a clean \belyi function, the edge segments come in pairs with a white vertex of valency equal to 2 between them. The union of a white vertex together with its pair of incident edge segments is called an {\em edge} and each edge corresponds to a unique disjoint 2-cycle of $g_1$. Thus, the black vertices and edges of $\mathcal G(\beta)$ are determined by the disjoint cycles of $g_0$ and $g_1$, respectively, and vice versa.

\begin{figure}[p]
\begin{center}
\begin{picture}(240,400)
\put(0,0){\epsfig{file=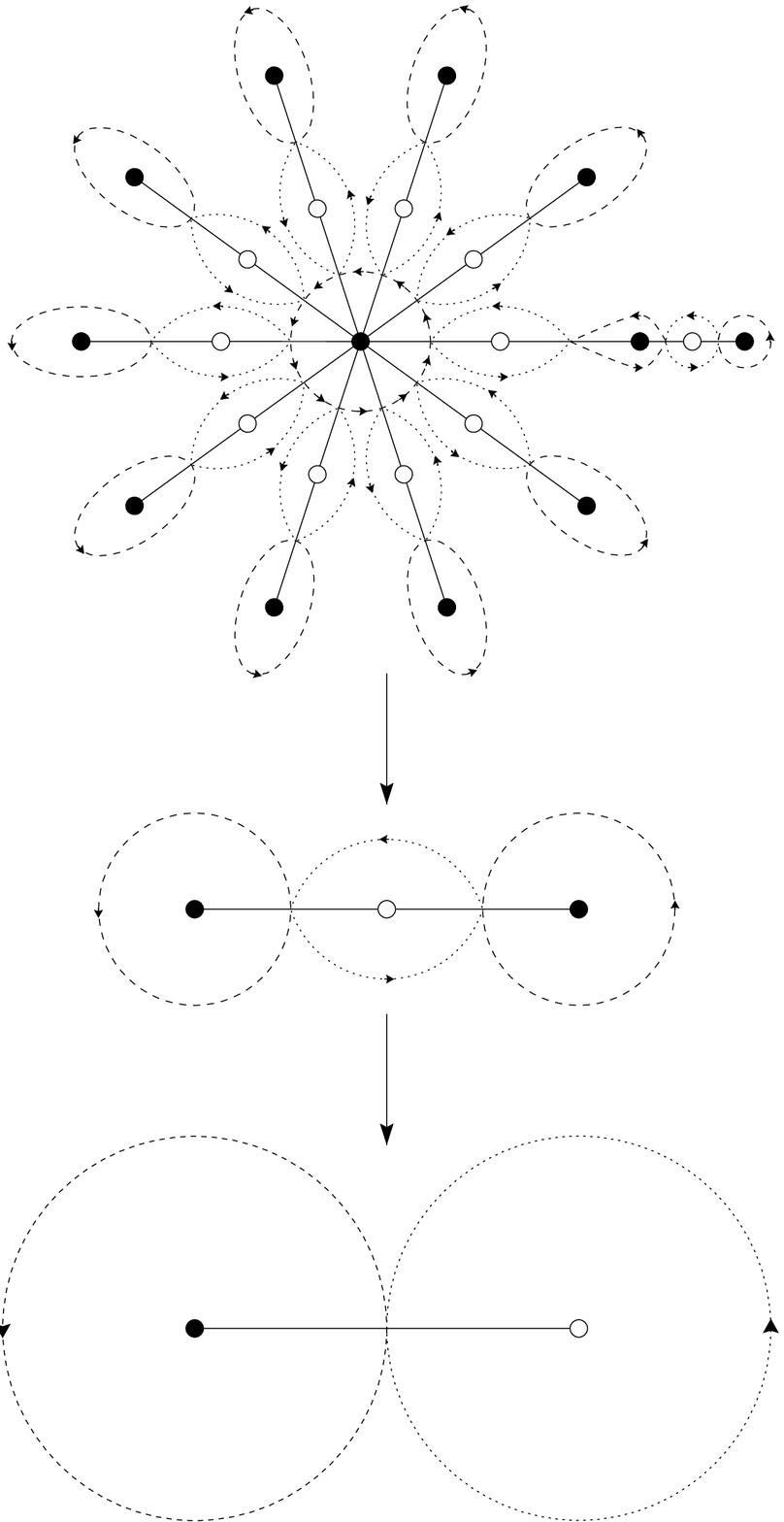,height=6in,width=3.5in}}

\put(-25,53){$\mathcal G_1$}
\put(57,45){$0$}
\put(115,46){$p=\frac{1}{2}$}
\put(190,45){$1$}
\put(130,125){$\beta_{1,1}$}

\put(-25,172){$\mathcal G(\beta_{1,1})$}
\put(57,165){$0$}
\put(85,169){$p_1$}
\put(112,166){$w=\frac{1}{2}$}
\put(160,169){$p_2$}
\put(190,165){$1$}
\put(130,223){$\beta_{10,1}$}

\put(-25,338){$\mathcal G(\psi)$}
\put(111,332){\footnotesize $0$}
\put(208,328){\footnotesize $\frac{10}{11}$}
\put(243,328){\footnotesize $1$}
\put(132,337){\footnotesize $\tilde p_1$}
\put(128,347){\footnotesize $\tilde p_2$}
\put(119,351){\footnotesize $\tilde p_3$}
\put(108,350){\footnotesize $\tilde p_4$}
\put(100,345){\footnotesize $\tilde p_5$}
\put(96,337){\footnotesize $\tilde p_6$}
\put(99,328){\footnotesize $\tilde p_7$}
\put(108,321){\footnotesize $\tilde p_8$}
\put(119,321){\footnotesize $\tilde p_9$}
\put(128,328){\footnotesize $\tilde p_{1\!0}$}

\put(183,342){$\tilde p_{11}$}
\put(175,363){$\tilde p_{12}$}
\put(145,388){$\tilde p_{13}$}
\put(79,391){$\tilde p_{14}$}
\put(50,363){$\tilde p_{15}$}
\put(47,328){$\tilde p_{16}$}
\put(175,307){$\tilde p_{20}$}
\put(144,283){$\tilde p_{19}$}
\put(80,283){$\tilde p_{18}$}
\put(50,308){$\tilde p_{17}$}

\put(213,348){$\tilde p_{21}$}
\put(232,348){$\tilde p_{22}$}

\end{picture}
\end{center}
\caption{The graphs $\mathcal G(\beta_{1,1})$ and  $\mathcal G(\psi)$ obtained by lifting the graph $\mathcal G_1$.}
\label{monodromy}
\end{figure}

Figure \ref{monodromy} displays the graphs $\mathcal G(\beta_{1,1})$ and $\mathcal G(\psi)$, where $\psi=\beta_{1,1} \circ \beta_{10,1}$, obtained by lifting the graph $\mathcal G_1$, which consists of a single edge $(0,1)$, a black vertex located at 0, and a white vertex located at 1. The graphs $\mathcal G(\beta_{1,1})$ and $\mathcal G(\psi)$ are shown in relation  to the lifts of rotations about 0 and 1 based at $p=\frac{1}{2}$, represented by the dashed and dotted loops, respectively. With the labelling of the points in $\Omega=\beta_{1,1}^{-1}(p)$ and $\tilde{\Omega}=\psi^{-1}(p)$ as in the figure, we have the following monodromy generators: 
\begin{align*}
G(\beta_{1,1}): g_0=&(1)(2),\\
g_1=&(1,2).\\
\\
G(\beta_{1,1} \circ \beta_{10,1}): \tilde{g}_0=&(1,2,3,4,5,6,7,8,9,10)(11,21),\\
\tilde{g}_1=&(1,11)(2,12)(3,13)(4,14)(5,15)(6,16)\\
&(7,17)(8,18)(9,19)(10,20)(21,22).\\
\end{align*}

\noindent The \belyi function $\beta_{1,1}$ is clean and has only one ramification point located at the white vertex $w=\frac{1}{2} \in \mathcal G(\beta_{1,1})$, corresponding to the only disjoint 2-cycle of $g_1$. The ramification order at $w$ is the cycle length $l_{1,1}=2$, which accounts for the two copies of $\mathcal G_1$ emanating from $w$. The \belyi function $\psi$ is also clean and $\mathcal G(\psi)$ has two copies of $\mathcal G_1$ emerging from every white vertex. In addition, $\psi$ has two ramification points at the black vertices located at 0 and $\frac{10}{11}$ in $\mathcal G(\psi)$, corresponding to the two non-trivial cycles of $\tilde{g}_0$. The black vertex located at 0 in $\mathcal G(\psi)$ corresponds to the 10-cycle in $\tilde{g}_0$ and has valency $\tilde{l}_{01}=10$. This accounts for the 10 copies of $\mathcal G(\beta_{1,1})$ emanating from 0, forming a floral pattern in $\mathcal G(\psi)$ centered at 0, which will be referred to as a 10-fold $\mathcal G(\beta_{1,1})$-bouquet. The black vertex at $\frac{10}{11} \in \mathcal G(\psi)$ corresponding to the 2-cycle of $\tilde{g}_0$, has valency $\tilde{l}_{02}=2$ and is the center of a 2-fold $\mathcal G(\beta_{1,1})$-bouquet.

\begin{figure}[p]
\begin{center}
\begin{picture}(200,516)
\put(0,0){\epsfig{file=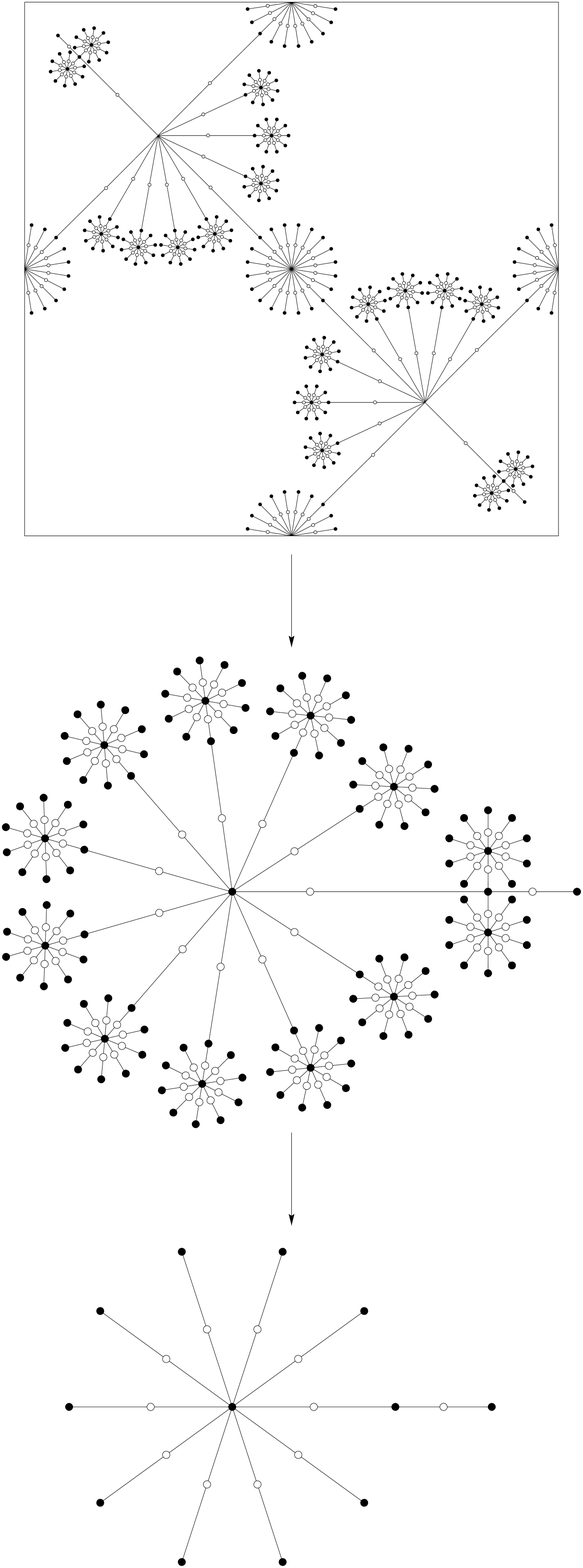,height=7.3in,width=3.0in}}

\put(-45,49){$\mathcal G(\psi)$}
\put(84,44){\footnotesize $0$}
\put(147,45){\footnotesize $\frac{10}{11}$}
\put(183,46){\footnotesize $1$}

\put(-45,223){$\mathcal G(\psi \circ f)$}
\put(81,221){\footnotesize $0$}
\put(183,222){\footnotesize $1$}
\put(215,219){\footnotesize $\frac{12}{11}$}

\put(-45,435){$\mathcal D(\beta^{2,7,11})$}
\put(113,324){$\pi^{2,7,11}$}
\put(113,129){$f$}
\end{picture}
\end{center}
\caption{The elliptic dessin $\mathcal D(\beta^{2,7,11})$ obtained by lifting the graph $\mathcal G(\psi \circ f)$ by $\pi^{2,7,11}$.}
\label{ellipticDessin}
\end{figure}

Figure \ref{ellipticDessin} shows the graphs $\mathcal G(\psi \circ f)$ and $\mathcal G(\beta^{2,7,11})$ obtained by lifting $\mathcal G(\psi)$ by $f$ and $f \circ \pi^{2,7,11}$, respectively. Unlike the previous figure, the dashed and dotted lines representing the lifts of rotations about 0 and 1 have been omitted, in order to emphasize the graphs; however, no monodromy information is lost by doing so, because $G(\beta^{2,7,11})$ is determined up to conjugacy by $\mathcal G(\beta^{2,7,11})$. The ramification point of $f$ at $v=1$ lies above the black vertex $f(1)=\frac{10}{11} \in \mathcal G(\psi)$, which joins a 10-fold $\mathcal G(\beta_{1,1})$-bouquet to the left and a 1-fold $\mathcal G(\beta_{1,1})$-bouquet to the right. Since $v$ has ramification order equal to 2, a lift of a neighbourhood of $\frac{10}{11} \in \mathcal G(\psi)$ by $f$ to a neighborhood of $v$ yields two 10-fold $\mathcal G(\beta_{1,1})$-bouquets and two 1-fold $\mathcal G(\beta_{1,1})$-bouquets emanating from $v$. This accounts for two copies of $\mathcal G(\psi)$ in $\mathcal G(\psi \circ f)$, which intersect at $v$; one starting at 0, which passes through $v$ and proceeds downward, and the other starting at $\frac{12}{11}$, which passes through $v$ and proceeds upward. The ramification point of $f$ at $0 \in \mathcal G(\psi \circ f)$ has ramification order equal to 11 and lies above $f(0)=1$, giving rise to the 11-fold $\mathcal G(\psi)$-bouquet centered at 0. The central black vertices of the twelve 10-fold $\mathcal G(\beta_{1,1})$-bouquets in $\mathcal G(\psi \circ f)$ are the zeros $\{ r_1, \ldots, r_{12}\}$ of $f$ and are ordered according to Figure \ref{fzeros}. 
The projection $\pi^{2,7,11}$ is ramified over only three points of $\mathcal G(\psi \circ f)$: $r_2, r_7$, and $r_{11}$, each having ramification order equal to 2. Hence, $\pi^{2,7,11}$ lifts each of the three 10-fold $\mathcal G(\beta_{1,1})$-bouquets centered at $r_2, r_7$, and $r_{11}$, to 20-fold $\mathcal G(\beta_{1,1})$-bouquets in $\mathcal G(\beta^{2,7,11})$. Since $\pi^{2,7,11}$ is a degree 2 covering and is unbranched away from $\{ r_2, r_7, r_{11} \}$, the remaining nine 10-fold $\mathcal G(\beta_{1,1})$-bouquets in $\mathcal G(\psi \circ f)$ are each lifted to two distinct 10-fold $\mathcal G(\beta_{1,1})$-bouquets giving a total of eighteen 10-fold $\mathcal G(\beta_{1,1})$-bouquets in $\mathcal G(\beta^{2,7,11})$. The dessin $\mathcal D(\beta^{2,7,11})$ consists of the graph $\mathcal G(\beta^{2,7,11})$ and a single face lying on the torus obtained by identifying the opposite sides of a square, as shown at the top of Figure \ref{ellipticDessin}.

%% file: GaloisActions.tex
\section{Galois Actions on Dessins}
$A_5$ is the orientation preserving symmetry group of
the regular dodecahedron and its action on the faces yields a permutation
representation $A_5 \rightarrow S_{12}.$ With the faces labelled as in
Figure 4, $A_5$ is given the following presentation:
\begin{equation*}
A_5 = \langle a,b \, |\, a^5=b^3=(ab)^2=1 \rangle,
\end{equation*}
where
 \begin{align*} 	
a &=
(2,3,4,5,6)(7,8,9,10,11),\\ 	
b &= (1,2,3)(4,6,7)(5,11,8)(9,10,12),\\
ab &= (1,2)(3,6)(4,11)(5,7)(8,10)(9,12).
\end{align*}

\begin{figure}[htb]
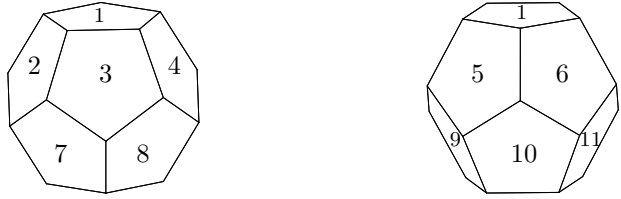

\label{dodecahedron}
\begin{picture}(200,75)
\put(76,65){\footnotesize $1$}
\put(51,45){$2$}
\put(78,42){$3$}
\put(104,45){$4$}
\put(61,12){$7$}
\put(92,12){$8$}
\begin{minipage}[b]{.46\linewidth}
\centering\epsfig{figure=12gon.ps,height=1in,width=1in}
\end{minipage}\hfill
\put(77,66){\footnotesize $1$}
\put(60,42){$5$}
\put(92,42){$6$}
\put(52,18){\footnotesize $9$}
\put(75,12){$10$}
\put(101,18){\footnotesize $11$}
\begin{minipage}[b]{.46\linewidth}
\centering\epsfig{figure=12gon2.ps,height=1in,width=1in}
\end{minipage}
\end{picture}
\caption{Two views of the dodecahedron are shown; the view on the right is a rotation of the view on the left by $\pi$ radians about the vertical axis. In both views, 12 is assigned to the hidden face at the bottom.}
\end{figure}

By taking the ordering of the roots $\{r_i\}_{i=1}^{12}$ of $f(x)$ in Figure \ref{fzeros}, and
identifying the $i^{th}$ root with the $i^{th}$ face of the dodecahedron, we get a Galois representation $A_5 \rightarrow G$.

Let $\mathcal{E}$ be the family of all elliptic curves $E_{i,j,k}$ as in (\ref{Eijk}), then $G = \mbox{Gal}(\mathbb Q,f)$ acts
on $\mathcal{E}$ as follows: an automorphism $\sigma \in G$ takes $E_{i,j,k}$ to a
conjugate curve $E_{i,j,k}^{\sigma}$ defined by
$$ E_{i,j,k}^{\sigma} = \left\{ [x,y,z]\in \mathbb{P}^2(\mathbb{C}) \, | \,
y^2=(x-\sigma(r_i)z)(x-\sigma(r_j)z)(x-\sigma(r_k)z) \right\}.  $$
By identifying the $i^{th}$ root $r_i$ with its index $i$, we have and induced action of $\sigma$ on the triple $(i,j,k)$, allowing us to rewrite $E_{i,j,k}^{\sigma}$ in the form
$$ E_{i,j,k}^{\sigma} = \left\{ [x,y,z]\in \mathbb{P}^2(\mathbb{C}) \, | \,
y^2=(x-r_{\sigma(i)}z)(x-r_{\sigma(j)}z)(x-r_{\sigma(k)}z) \right\}. $$
So the action of $\sigma$ on the triple $(i,j,k)$ induces an action on the elliptic curve $E_{i,j,k}$, and further, induces an action on the dessin $\mathcal D(\beta^{i,j,k})$ associated to the \belyi pair $(E_{i,j,k},\beta^{i,j,k})$. For example, $a : E_{2,7,11} \mapsto E_{3,8,7} = E_{3,7,8}$, which induces the action $a : D(\beta^{2,7,11}) \mapsto D(\beta^{3,7,8})$. Moreover, the orbits containing the triple $(2,7,11)$ under the action of the subgroups $<\!\!a\!\!>$, $<\!\!b\!\!>$, and $<\!\!ab\!\!> \subset A_5$ are as follows:
\begin{align*} 	
<\!\!a\!\!> &: \; \{ (2,7,11),(3,7,8),(4,8,9),(5,9,10),(6,10,11) \},\\
<\!\!b\!\!> &: \; \{ (2,7,11),(3,4,8),(1,5,6) \},\\
<\!\!ab\!\!> &: \; \{ (2,7,11),(1,4,5) \}.
\end{align*}
The orbits of the corresponding dessins are illustrated in Figures 5, 6, and 7.

\begin{figure}[p]
\label{actiona}
\begin{center}
\begin{picture}(345,515)
\put(0,0){\epsfig{file=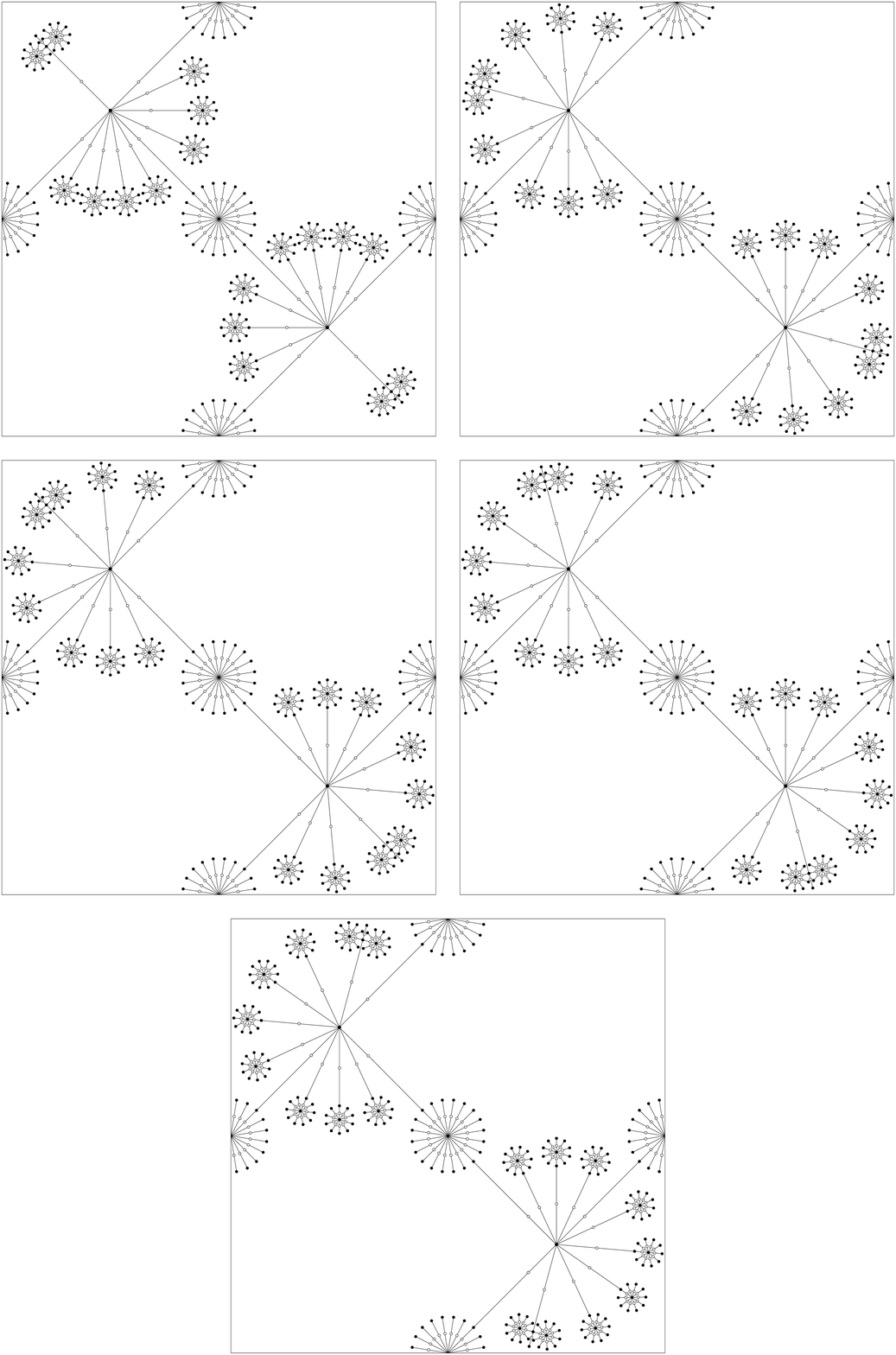,height=7.3in,width=4.9in}}
\end{picture}
\end{center}
\caption{The orbit of the dessin $\mathcal D(\beta^{2,7,11})$ under the action of $<\!\!a\!\!> \subset G$, shown from left to right as successive powers of $a$ act on $\mathcal D(\beta^{2,7,11})$. }
\end{figure}

\begin{figure}[p]
\label{actionb}
\begin{center}
\begin{picture}(345,280)
\put(0,0){\epsfig{file=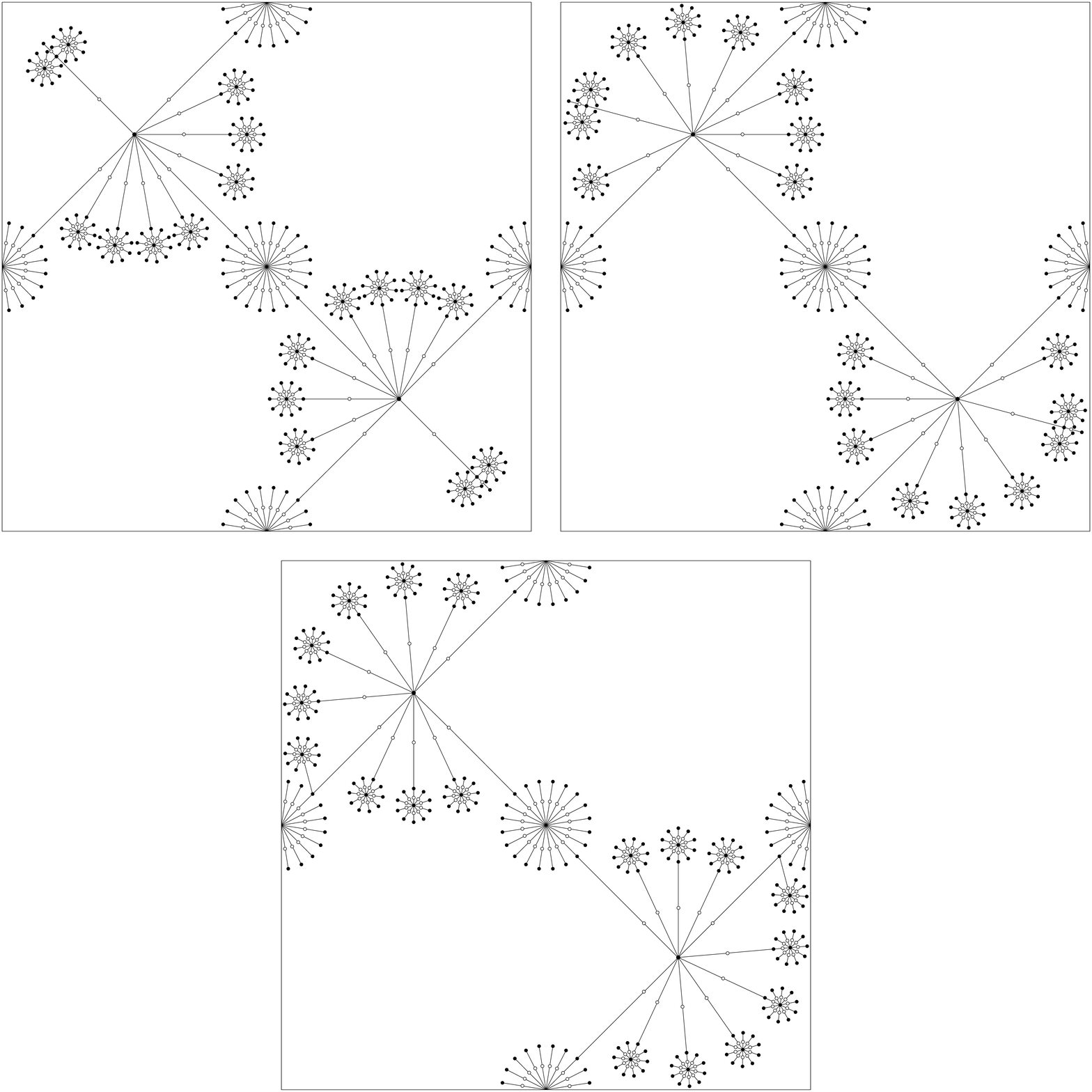,height=4.5in,width=4.9in}}
\end{picture}
\end{center}
\caption{The orbit of the dessin $\mathcal D(\beta^{2,7,11})$ under the action of $<\!\!b\!\!> \subset G$, shown from left to right as successive powers of $b$ act on $\mathcal D(\beta^{2,7,11})$. }
\end{figure}

\begin{figure}[p]
\label{actionab}
\begin{center}
\begin{picture}(345,180)
\put(0,0){\epsfig{file=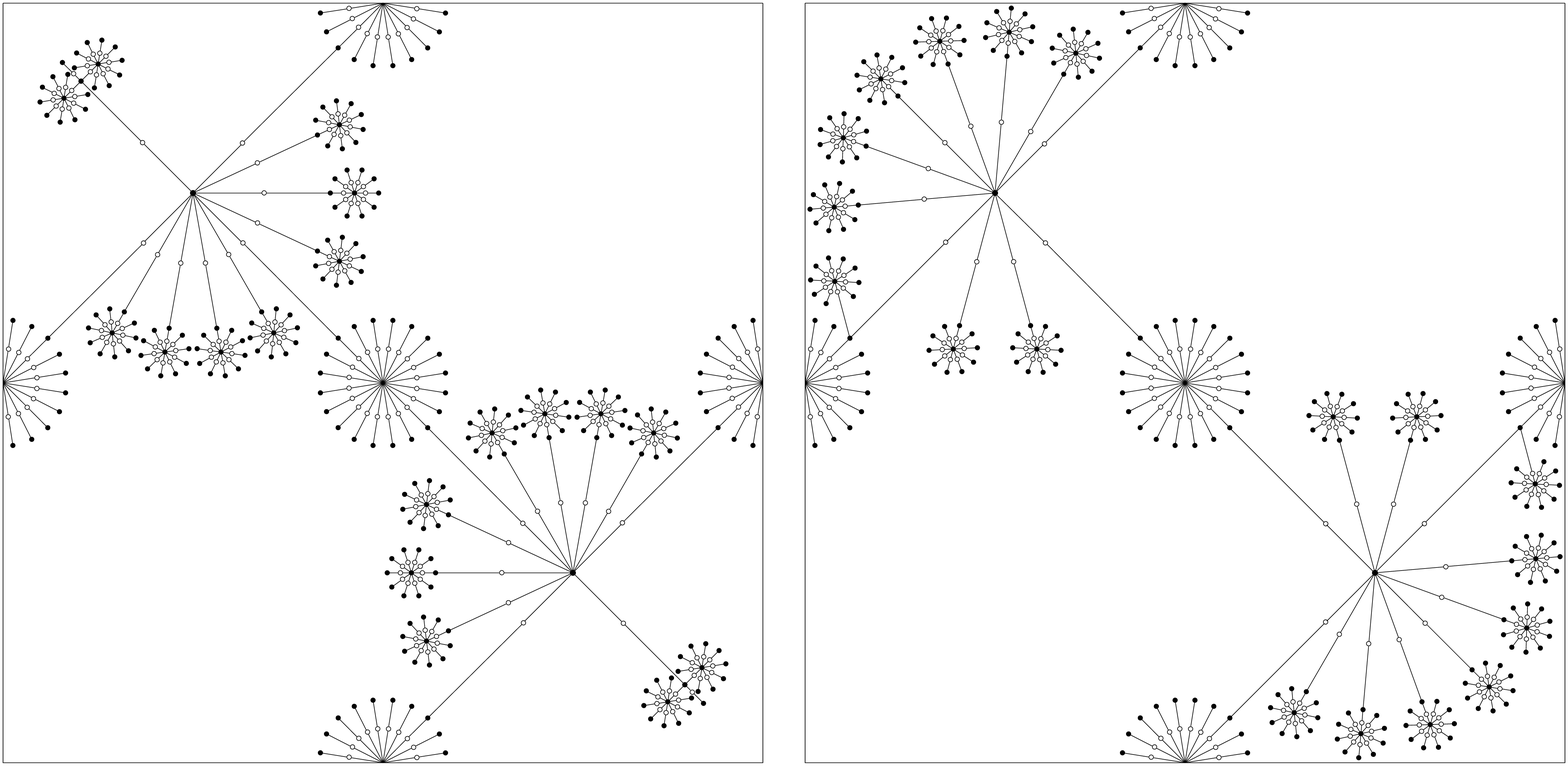,height=2.3in,width=4.9in}}
\end{picture}
\end{center}
\caption{The action of $ab$ on the dessin $\mathcal D(\beta^{2,7,11})$.}
\end{figure}